\documentclass[reqno,a4paper,11pt]{amsart}
\usepackage[all,poly]{xy}
\usepackage{amsfonts}
\usepackage[mathcal]{eucal}
\usepackage{amssymb}
\usepackage{amsmath}
\usepackage{mathrsfs}
\usepackage{color}
\usepackage[pagebackref,colorlinks]{hyperref}
\usepackage{enumerate}
\usepackage{comment}
\usepackage{hhline}
\usepackage{graphics}
\usepackage{float}
\usepackage{comment}

\usepackage{enumerate}

\theoremstyle{plain}

\theoremstyle{definition}

\newcommand{\Aut}{\operatorname{Aut}}

\newcommand{\chr}{\operatorname{char}}

\title{A note on multiplicative commutators of division rings}

\author{Roozbeh Hazrat}
\address{Western Sydney University,
Australia}\email{r.hazrat@westernsydney.edu.au}

\thanks{The author acknowledges Australian Research Council grant DP160101481. This work was done at the University of M\"unster, where the author was a Humboldt Fellow.}

\subjclass[2000]{} 

\begin{document}
\maketitle

\begin{abstract} 
We give an example of a division ring $D$ whose multiplicative commutator subgroup does not generate $D$ as a vector space over its centre, thus disproving the conjecture posed in \cite{aa}. 
\end{abstract}

Let $D$ be a division ring. Denote by $F:=Z(D)$, the centre of $D$ and $D'$ its multiplicative commutator subgroup, i.e., the group generated by the set of multiplicative commutators $\big\{xyx^{-1}y^{-1} \mid x,y \in D\backslash \{0\} \big\}$. There are classical results due to Herstein, Kaplansky and Scott, among others, showing that the group $D'$ is ``dense'' in $D$~(see for example \cite[\S13]{lam}). In \cite{aa} the authors study the $F$-vector space $T(D)$ generated by the set of multiplicative commutators.
They prove that if $T(D)$ is radical over $F$, then $D = F$, and if $\dim_F T(D) < \infty$, then $\dim_F D < \infty$. Furthermore, they prove that $T(D)$ contains all separable elements of $D$ and thus if $D$ is an algebraic division ring over its centre with $\chr(D) = 0$, then  $T(D)=D$.

They then conjecture \cite[Abstract and Conjecture~1]{aa} that a division ring is generated by all multiplicative commutators as a vector space over its centre, i.e., $T(D)=D$ for any arbitrary division ring. 

Here we give a counterexample to this conjecture. In fact we show that the multiplicative subgroup $D'$ can not recover $D$ as a vector space and $[D:T(D)]=\infty$.

We recall the Hilbert classical construction of division rings (see \cite[\S1]{draxl}). Let $L$ be a field, $\sigma \in \Aut(L)$ and
$F$ be the fixed field of $\sigma$. Let $D=L((x,\sigma))$ be the division ring of formal Laurent series, consisting of elements  
$\sum_{i=n}^{\infty}a_it^i$, where $a_i\in L$ and $n\in \mathbb Z$, addition defined component-wise and multiplication by
\begin{equation}\label{rh1}
\Big(\sum_{i=m}^{\infty}a_it^i\Big)\Big(\sum_{j=n}^{\infty}b_jt^j\Big)=\sum_{r=n+m}^{\infty}\Big(\sum_{i+j=r}a_i\sigma^{i}(b_j)\Big)t^r.
\end{equation}

By~\cite[\S1, Lemma 4]{draxl} if the order of $\sigma$ is infinite then  $Z(D)=F$ and $[D:F]=\infty$ and if the order is a finite number $n$, then $Z(D)=F((t^n))$ and $[D:Z(D)]=n^2$.  It is easy to see that $D=L((x,\sigma))$ is a valued division ring with value group $\mathbb Z$ as follows  \begin{align*}
v:D^* &\longrightarrow \mathbb Z,\\
\sum_{i=n}^{\infty}a_it^i &\longmapsto n.
\end{align*}
Therefore we have 
\begin{equation}\label{rh2}
D' \subseteq \ker (v)= \Big \{L +\sum_{i>0}a_it^i \Big\} .
\end{equation}
Now if we choose an automorphism  $\sigma$ of infinite order, since $Z(D)=F$, Equation~\ref{rh2} shows that $T(D) \subseteq \Big \{L + \sum_{i>0}a_it^i \Big\}$, and thus $[D:T(D)]=\infty$.

\end{document}